\documentclass[11pt]{amsart}
\usepackage[all,cmtip]{xy}
\usepackage{amsmath,amssymb,amscd,mathrsfs,graphicx}
\usepackage[hidelinks]{hyperref}
\newtheorem{theorem}{Theorem}
\newtheorem{lemma}{Lemma}
\newtheorem{proposition}{Proposition}
\newtheorem{corollary}{Corollary}
\theoremstyle{definition}
\newtheorem{definition}{Definition}

\theoremstyle{remark}
\newtheorem{remark}{Remark}
\DeclareMathOperator{\Nef}{Nef}
\DeclareMathOperator{\Mov}{Mov}
\DeclareMathOperator{\Eff}{Eff}
\DeclareMathOperator{\Bs}{Bs}
\DeclareMathOperator{\BF}{BF}
\DeclareMathOperator{\Coker}{Coker}
\DeclareMathOperator{\Pic}{Pic}
\DeclareMathOperator{\Cl}{Cl}
\DeclareMathOperator{\Cart}{Cart^{B}}
\newcommand{\MD}{\mathcal D}
\newcommand{\RG}{\mathop{\mathcal R^{G}}\nolimits}
\newcommand{\RGX}[1]{\mathcal R^{G}_{#1}}
\DeclareMathOperator{\MP}{\mathcal P}

\begin{document}
\title{Gale Duality and the Movable Cone of Spherical Varieties}
\author{V. S. Zhgoon}
\address{Moscow Institute of Physics and Technology (National Research University), Dolgoprudny, Moscow Region; HSE University, Moscow; Federal Research Center Scientific Research Institute for System Analysis of the National Research Centre Kurchatov Institute, Moscow}
\email{zhgoon@mail.ru}
\date{September 5, 2026}
\thanks{This research was supported by the Russian Science Foundation (project 25-21-00358).}
\begin{abstract}
We describe the movable cone  of a spherical variety and its decomposition into chambers in terms of Gale duality.
\end{abstract}
\maketitle
\section{Introduction}
Gale duality arises naturally in problems of toric geometry and topology (see \cite{Arzh,panov,Sh}). In \cite{Gale}, it was used to describe the nef cone of a spherical variety.

In the present paper, this approach is applied to the study of movable divisors on spherical varieties. We describe the monoid and the cone of movable divisors in terms of Gale duality, generalizing the known description in the toric case. We establish a connection between the mutual position of Gale dual cones and the combinatorics of colored fans. On this basis, we describe the decomposition of the movable cone into chambers that are the nef cones of various equivariant birational models isomorphic in codimension one.

The foundations of Mori theory for spherical varieties are presented in Brion's work \cite{BrionMori}. For the basic notions of the theory of spherical varieties, we refer the reader to \cite{Brion,Knop,Tim}.

The author is sincerely grateful to D. A. Timashev, V. L. Popov, T. E. Panov, V. V. Shokurov, and K. V. Loginov for useful discussions.

Let $G$ be a reductive group over an algebraically closed field $k$ (to avoid unnecessary qualifications, we assume that $k$ has characteristic zero). Let $B$ be a Borel subgroup of $G$, and let $X$ be a normal spherical variety, that is, an irreducible normal algebraic variety with an action of $G$ containing an open $B$-orbit. It is well known that $X$ contains only finitely many $G$-orbits and even finitely many $B$-orbits (see \cite{Vin,Tim}).

\section{Movable divisors}
\begin{definition}
A Weil divisor $D$ is called \emph{movable} if it is effective and the base locus of the complete linear system of some positive integral multiple of $D$ has codimension at least $2$.
\end{definition}
For a Weil divisor $D$, the complete linear system $|D|$ consists of effective divisors linearly equivalent to $D$, and its base locus is the intersection of their supports. On the smooth locus of $X$, this is the usual base locus of a linear system.
We denote by $\Mov(X)$ the rational cone in $\Cl(X)_{\mathbb Q}$ generated by the classes of movable divisors.
Let $\MP_{\Mov(X)}$ denote the set of classes of effective Weil divisors whose complete linear systems have base loci of codimension at least $2$.
The inclusion $|D_1|+|D_2|\subseteq |D_1+D_2|$ implies that the base locus of the sum of two divisors is contained in the union of their base loci; thus $\MP_{\Mov(X)}$ is a monoid.

A Cartier divisor is called \emph{semiample} if some positive integral multiple of it is base-point-free.

For a spherical variety $X$, it will be convenient to use another definition. We precede its formulation with the following proposition.

\begin{proposition}\label{prop:movdef}
Let $D$ be a movable effective Cartier divisor on a complete normal variety $X$ equipped with an action of a connected linear algebraic group $G$. Then there exist an integer $n>0$, a proper $G$-equivariant birational morphism $\psi :\widetilde{X} \rightarrow X$ from a normal variety, and an effective nef Cartier divisor $\widetilde{D}$ on $\widetilde{X}$ such that $nD=\psi_*\widetilde{D}$ (where the pushforward is understood in the sense of Weil divisors). Moreover, $\psi^*(nD)=\widetilde{D}+F$, where $F$ is an effective exceptional divisor whose image $\psi(\operatorname{Supp} F)$ is contained in $\Bs(|nD|)$ ($\psi_*(F)=0$).

Conversely, if $X$ is spherical, then every Cartier divisor $D$ such that $nD=\psi_*\widetilde{D}$ for some $n>0$ and some proper $G$-equivariant birational morphism $\psi:\widetilde X\to X$ from a normal variety, where $\widetilde D$ is an effective nef Cartier divisor, is movable.
\end{proposition}
\begin{proof}
After passing to a covering $\widetilde{G}\rightarrow G$, we may assume that the invertible sheaf $\mathcal O(D)$ is $G$-linearized. Choose $n>0$ so that $\Bs(|nD|)$ has codimension at least $2$. The base ideal $\mathfrak b$ of the complete linear system $|nD|$ is $G$-invariant.

Let $\psi:\widetilde X\to X$ be the normalization of the blowup of the ideal $\mathfrak b$. Then $\psi$ is a proper $G$-equivariant birational morphism that is an isomorphism over $X\setminus \Bs(|nD|)$, and
\[\mathfrak b\mathcal O_{\widetilde X}=\mathcal O_{\widetilde X}(-F),\]
where $F$ is an effective Cartier divisor. Since $\Bs(|nD|)$ has codimension at least $2$, the divisor $F$ is exceptional and $\psi_*(F)=0$.

The complete linear system $|nD|$ defines a map
\[
f:X\setminus \Bs(|nD|)\longrightarrow \mathbb P(H^0(X,\mathcal O(nD))^*).
\]
Let $Y$ be the closure of the image of $f$. The morphism $\psi$ resolves the indeterminacy of this map, yielding the commutative diagram
\[
\xymatrix{
   &  &  \widetilde{X} \ar[dl]_{\psi} \ar[rd]^{\widetilde f} & \\
   & X  &\ar@{_{(}->}[l] X\setminus \Bs(|nD|) \ar[r]^{\  \  \  \  \  f} & Y. }
\]
Here $\widetilde f:\widetilde X\to Y$ is a morphism extending $f\circ\psi$.

Set $\widetilde D=\psi^*(nD)-F$. By the definition of the base ideal, the pulled-back sections of $\mathcal O(nD)$ generate the sheaf $\mathcal O_{\widetilde X}(\widetilde D)$ after removal of the common fixed part $F$. Consequently, $\widetilde D$ is base-point-free and nef. The section defining the effective divisor $nD$ defines the effective divisor $\widetilde D$ after removal of $F$. Thus $\psi^*(nD)=\widetilde D+F$ and $\psi_*\widetilde D=nD$.

We now prove the converse. Since $\widetilde X$ is also spherical, by \cite[Theorem~2]{Gale} (see also \cite{Brion,BrionPic,Tim}) the nef Cartier divisor $\widetilde D$ is base-point-free. By normality of $X$ and properness and birationality of $\psi$, there is an open subset $U\subset X$ whose complement has codimension at least $2$, over which $\psi$ is an isomorphism. On $U$, the divisors $\widetilde D$ and $nD$ are identified. The sections of $\mathcal O_{\widetilde X}(\widetilde D)$ restricted to $\psi^{-1}(U)$ extend to sections of $\mathcal O_X(nD)$ by normality of $X$. Hence $|nD|$ is base-point-free on $U$, proving that $D$ is movable.
\end{proof}
\begin{remark}
The assertion of Proposition~\ref{prop:movdef}, except for the equality $\psi^*(nD)=\widetilde D+F$, remains valid for Weil divisors. In this case, a linearization of $\mathcal O_X(D)$ is obtained by extending the linearization of its restriction to the smooth locus of $X$ (after passing, if necessary, to a covering of $G$).

For $\widetilde X$ one takes the normalization of the closure of the graph of the rational map defined by the complete linear system $|nD|$, and for $\widetilde D$ the pullback of the hyperplane section corresponding to the effective divisor $nD$. Since $nD$ and $\widetilde D$ are identified over an open subset whose complement has codimension at least two, we have $nD=\psi_*\widetilde D$. The proof of the converse remains unchanged.
\end{remark}
By Proposition~\ref{prop:movdef}, we may use the following definition of a movable divisor.
\begin{definition}\label{def:mov}
Let $X$ be a complete spherical variety. A Weil divisor $D$ is called movable if there exist a proper $G$-equivariant birational morphism $\psi :\widetilde{X} \rightarrow X$ and an effective nef Cartier divisor $\widetilde{D}$ on a normal variety $\widetilde{X}$ such that $nD=\psi_*\widetilde{D}$ for some integer $n>0$.
\end{definition}

\section{Colored fans and Gale duality}

We recall the necessary notation from the theory of spherical varieties.
Let $k(X)$ denote the field of rational functions on $X$. The existence of an open $B$-orbit implies that the field of rational $B$-invariant functions is trivial: $k(X)^B=k$.

Let $\Xi(B)=\operatorname{Hom}(B,k^\times)$ be the character lattice of $B$. Consider the weight lattice $\Xi(X)$ of $X$, that is, the set of weights of rational $B$-semi-invariant functions on $X$:
\[
\Xi(X)= \{\chi \in \Xi(B) \mid \exists\, f_\chi\in k(X)\setminus \{0\},\
f_\chi(b^{-1}x)=\chi(b) f_{\chi}(x) \}.
\]
We denote the dual lattice by $\Lambda(X)=\operatorname{Hom}(\Xi(X),\mathbb Z)$.
Let $\mathcal V_X$ denote the set of rational divisorial $G$-invariant valuations of $k(X)$ (the valuation cone).
Let $f:X\rightarrow Y$ be a dominant $G$-equivariant morphism of spherical $G$-varieties. It induces a homomorphism of character lattices $f^*:\Xi(Y)\rightarrow \Xi(X)$ arising from the map $f^*:k(Y)\rightarrow k(X)$, and the dual map $f_*:\Lambda(X) \rightarrow \Lambda(Y)$.

For a $G$-orbit $Y$, let $\RGX{X,Y}$ denote the set of valuations corresponding to $G$-invariant prime divisors containing $Y$. Let $\mathcal D$ denote the set of $B$-invariant prime divisors that are not $G$-invariant; such divisors are called \emph{colors}. Let $\mathcal D_Y$ denote the set of colors containing $Y$. We write $\nu_D$ for the valuation of $k(X)$ corresponding to a prime divisor $D$. Let $\RGX{X}$ denote the set of $G$-invariant valuations corresponding to $G$-invariant divisors on $X$, and let $\RGX{X,Y}$ be the subset of $\RGX{X}$ consisting of valuations corresponding to $G$-invariant divisors containing $Y$.

Consider the variety $X_Y$ consisting of the $G$-orbits whose closures contain $Y$. Since $X$ is normal, by Sumihiro's theorem (\cite{KKLV}) it admits an open covering by $G$-invariant quasiprojective varieties, each of which embeds into the projectivization of a vector space (possibly different for each member of the covering) with a linear action of $G$. Since $X_Y$ is the smallest open $G$-subvariety containing $Y$, it is contained in a member of this covering and is therefore quasiprojective.

We now turn to the maps used to define Gale duality. Consider the free abelian group $\mathbb Z\langle e_\nu,e_D\rangle$ generated by basis elements $e_\nu$ indexed by valuations $\nu\in\RGX{X}$ and basis elements $e_D$ indexed by divisors $D\in\mathcal D$. Consider the natural map
\[
\mathbb Z\langle e_\nu,e_D\rangle \xrightarrow{A} \Lambda(X),
\]
which assigns to $e_\nu$ and $e_D$ linear functionals on the character lattice $\Xi(X)$ whose values at $\chi\in\Xi(X)$ are given by
\[
A(e_{\nu})(\chi):= \nu(f_\chi), \qquad A(e_D)(\chi):= \nu_D(f_\chi),
\]
where $f_\chi$ is a rational $B$-semi-invariant function of weight $\chi$, determined up to multiplication by a constant since $k(X)^B=k$. The map $A$ embeds the set $\RGX{X}$ and maps the set $\mathcal D$ into the lattice $\Lambda(X)$.
Let $\mathcal C_{X,Y}$ denote the cone generated in the lattice $\Lambda(X)$ by the images $Ae_\nu$ for $\nu\in\RGX{X,Y}$ and $Ae_D$ for $D\in\mathcal D_Y$ (this cone is colored by the colors in $\mathcal D_Y$). If the closure of a $G$-orbit $Y_1$ contains $Y$, then $\mathcal C_{X,Y_1}$ is a face of $\mathcal C_{X,Y}$.

For any convex cone or polyhedron $P$, we denote its relative interior by $P^\circ$. We write $\langle v_1,\ldots,v_s\rangle_{\geqslant 0}$ for the convex cone generated by the vectors $v_1,\ldots,v_s$, and $\langle v_1,\ldots,v_s\rangle_{>0}$ for its relative interior, that is, the set of linear combinations of these vectors with strictly positive coefficients.

\begin{definition}
A face $\mathcal C$ of the cone $\mathcal C_{X,Y}$ is called \emph{colored} if its relative interior $\mathcal C^\circ$ has nonempty intersection with the valuation cone $\mathcal V_X$.
\end{definition}

In this case the face corresponds to a simple variety $X_{Y_1}$. The separatedness of $X$ implies that the colored cones $\mathcal C_{X,Y}$ form a so-called colored fan. Namely, for the relative interiors of any two distinct colored faces $\mathcal C_{X,Y_i}$, $i=1,2$, we have
\begin{equation}\label{eq:separ}
\mathcal C_{X,Y_1}^\circ \cap \mathcal C_{X,Y_2}^\circ\cap \mathcal V_X=\emptyset.
\end{equation}
It follows directly from the definition that the images of the elements $\nu$ corresponding to $G$-invariant divisors containing $Y$ form a set of colored rays of the colored cone $\mathcal C_{X,Y}$. We denote the colored fan by $\Sigma$, and the set of its rays with their corresponding markings by $\Sigma(1)$.

\begin{remark}
Note that a ray in $\Sigma(1)$, being colored, may be marked by a color whose valuation determines an element of $\Lambda(X)$ lying in the valuation cone. In this case, however, no $G$-invariant divisor on $X$ corresponds to this valuation.
\end{remark}

Consider also the dual map $A^*$:
\[
\Xi(X) \xrightarrow{A^*} \mathbb Z\langle X_\nu, D\rangle,
\]
\[
A^*(\chi) = (f_\chi)=\sum_{\nu \in \RGX{X}} n_\nu X_\nu+\sum_{D \in \mathcal D} n_D D.
\]
This map assigns to $\chi$ the collection of values $n_\nu=\nu(f_\chi)$ for all $\nu\in\RGX{X}$ and $n_D=\nu_D(f_\chi)$ for all $D\in\mathcal D$. The group $\mathbb Z\langle X_\nu,D\rangle$ contains the subgroup $\Cart(X)$ of $B$-invariant Cartier divisors. On the other hand, $\Cl(X)$ is generated by the finite set of divisors in the complement of the open $B$-orbit (whose class group is trivial). Consequently, the subgroup $\Pic(X)\subset\Cl(X)$ is also finitely generated. Passing to a covering $\widetilde G\rightarrow G$, we may assume that each generator of $\Cl(X)$ is represented by a $G$-linearized reflexive sheaf of rank one, as in the proof of Proposition~\ref{prop:movdef}.
Consider the canonical map $\pi$ from $\mathbb Z\langle X_\nu,D\rangle$ to the quotient group $\Cl(X)=\Coker A^*=\mathbb Z\langle X_\nu,D\rangle/A^*(\Xi(X))$. The image of $\Cart(X)$ under $\pi$ is $\Pic(X)$. We also denote by $\mathcal P_{\Eff}$ the monoid of classes of effective Weil divisors (and by $\Eff(X)$ the rational cone generated by this monoid). Since a linear system contains a $B$-semi-invariant section by the Lie--Kolchin theorem, every element of $\mathcal P_{\Eff}$ is a nonnegative linear combination of classes of $B$-invariant divisors.

\begin{definition}
The \emph{Gale dual monoid} of a cone $\mathcal C_{X,Y}$ of the colored fan $\Sigma_X$ is the monoid $\MP_{X,\widehat{Y}}\subset\operatorname{Coker} A^*$ generated by the nonnegative integral linear combinations of the images under $\pi$ of the rays $X_\nu$, where $\nu\in\RGX{X}\setminus\RGX{X,Y}$, and $D\in\mathcal D\setminus\mathcal D_Y$.
\end{definition}
Let $\Lambda_{X,\widehat{Y}}$ denote the lattice in $\operatorname{Coker} A^*$ generated by the images under $\pi$ of the elements $X_\nu$, where $\nu\in\RGX{X}\setminus\RGX{X,Y}$, and $D\in\mathcal D\setminus\mathcal D_Y$ (that is, the lattice generated by the points of the monoid $\MP_{X,\widehat{Y}}$). We denote by $\MP^{\mathbb Q}_{X,\widehat{Y}}$ the rational cone generated by the monoid $\MP_{X,\widehat{Y}}$.

The following theorem describes the Picard group and the nef cone.

\begin{theorem}\label{thm:pic}\label{Nef}
(\cite[Theorems~1 and 2]{Gale}) Let $X$ be a spherical variety over a field $\Bbbk$ of characteristic zero, and let $E$ be a Cartier divisor on $X$ invariant under $B$.
\begin{itemize}
\item[(i)] We have
\[
\Pic(X)=\pi(\Cart(X))=\bigcap_{Y\subset X}\Lambda_{X,\widehat{Y}},
\]
where $Y$ runs over the closed $G$-orbits in $X$.
\item[(ii)] The complete linear system $|E|$ is base-point-free on a $G$-orbit $Y$ if and only if the class of $E$ belongs to $\MP_{X,\widehat{Y}}$.
\item[(iii)] The complete linear system $|E|$ is base-point-free if and only if the class of $E$ belongs to the monoid
\[
\MP_{\BF}:=\bigcap_{Y\subset X}\MP_{X,\widehat{Y}},
\]
where the intersection is taken over all closed $G$-orbits in $X$. If $X$ is complete, this monoid coincides with the monoid of classes of nef Cartier divisors (that is, divisors having nonnegative intersection with effective complete curves).
\item[(iv)] The divisor $E$ is ample if and only if its class belongs to the interior of the rational cone generated by the monoid $\MP_{\BF}$ in the space $\Pic(X)_{\mathbb Q}$.
\end{itemize}
\end{theorem}

\begin{remark}
An interior point of a cone is a point belonging to the cone together with a neighborhood in $\Pic(X)\otimes_{\mathbb Z}\mathbb Q$.
\end{remark}

\section{Equivariant Morphisms and Subdivisions}

Recall how $G$-equivariant morphisms between spherical varieties are described in terms of colored fans.

\begin{proposition}[\cite{Knop}]
Let $f:\widetilde{X}\rightarrow X$ be a $G$-equivariant morphism of spherical $G$-varieties. Denote by $\mathcal D_f$ the set of colors of $\widetilde X$ that map dominantly onto $X$. The map $f_*:\Lambda(\widetilde X)\rightarrow\Lambda(X)$ can be applied to the colored cones $\mathcal C_{\widetilde X,\widetilde Y}$ of the fan of $\widetilde X$ (where $\widetilde Y\subset\widetilde X$ is a $G$-orbit). For a $G$-orbit $Y\subset X$ contained in the closure of $f(\widetilde Y)$, the following conditions must hold:
\begin{itemize}
\item[(i)] $f_*(\mathcal C_{\widetilde X,\widetilde Y})\subset\mathcal C_{X,Y}$,
\item[(ii)] $f_*(\mathcal D_{\widetilde Y}\setminus\mathcal D_f)\subset\mathcal D_Y$.
\end{itemize}
\end{proposition}

\begin{remark}
Since the restriction of the $G$-equivariant map $f$ to the open $G$-orbit in $\widetilde X$ is equidimensional and every color $D\in\mathcal D_{\widetilde X}$ has a nonempty open intersection with this orbit, either $D$ maps dominantly onto $X$, or the intersection of $f(D)$ with the open orbit in $X$ is a proper subset of codimension one in $X$ and hence determines a color. In the situation under consideration, where $f$ is birational, we have $\Xi(X)\cong\Xi(\widetilde X)$. Moreover, in this case the restriction of $f$ to the open $G$-orbit is an isomorphism, so the set $\mathcal D_f$ is empty.
\end{remark}

\begin{proposition}[\cite{Knop}]
A $G$-equivariant morphism $f:\widetilde X\rightarrow X$ of spherical $G$-varieties is proper if and only if
\[
\bigcup_{\widetilde{\mathcal C}\in\Sigma_{\widetilde X}}\widetilde{\mathcal C}\cap\mathcal V_{\widetilde X}
\supset
\bigcup_{\mathcal C\in\Sigma_X}f_*^{-1}(\mathcal C\cap\mathcal V_X).
\]
\end{proposition}

We have the following commutative diagram:
\[
\xymatrix{
   & \Xi(X) \ar[r] \ar[d]^{\cong} & \mathbb Z\langle X_\nu,D\rangle_{\nu \in \RG \, D \in \mathcal D_X} \ar[r] & \Cl(X) \\
   & \Xi(\widetilde{X}) \ar[r] & \mathbb Z\langle X_{\widetilde{\nu}},D\rangle_{\widetilde{\nu} \in \widetilde{\RG} \, D \in \mathcal D_{\widetilde{X}}} \ar[u]_{f_*} \ar[r] & \Cl(\widetilde{X}) \ar[u]_{f_*} },
\]
where the middle and right vertical arrows $f_*$ are quotient maps modulo the classes of the exceptional divisors of $f$. On the one hand, exceptional divisors are characterized by the fact that their images have higher codimension; on the other hand, they can be characterized as the colored rays in $\widetilde{\RG}$ corresponding to $G$-invariant divisors and not contained in $\RG$.

\begin{lemma}\label{baricenter}
For a colored fan $\Sigma$ of a spherical variety $X$, there exists a subdivision into colored cones such that every colored cone of maximal dimension in the subdivision contains at most one colored ray of the original fan.
\end{lemma}

\begin{proof}
We construct the subdivision starting with the cones of highest dimension. Namely, for each colored cone $\mathcal C_0$ of maximal dimension, choose a vector $v_0$ in the intersection of its relative interior with the valuation cone $\mathcal V_X$. Label it by a $G$-invariant divisor (even if there exists a color with the same label). By convexity, the cone is the union of the cones generated by $v_0$ and a facet of $\mathcal C_0$.

Next, applying the same procedure to a colored facet $\mathcal C_1$, we obtain a vector $v_1$ and a subdivision of $\mathcal C_1$ into cones having the common ray generated by $v_1$. We then replace the cone generated by $v_0$ and $\mathcal C_1$ by the union of the cones generated by facets of $\mathcal C_1$ and the vectors $v_0$ and $v_1$. As a result of this inductive procedure, in the relative interior of each colored face $\mathcal V_i$ we find a vector $v_i\in\mathcal C_X$. We then replace $\mathcal C_i$ by the union of cones, each generated by a facet $\mathcal C_{i+1}$ of $\mathcal C_i$ and the vector $v_i$. Thus we add to the subdivision of the fan cones of maximal dimension generated by $v_0,v_1,\ldots,v_i$ and the face $\mathcal C_{i+1}$. If the face $\mathcal C_{i+1}$ is colored, these cones are also colored.

This process terminates when we reach a one-dimensional colored face. If $X_\nu\in\mathcal V_X$ is a ray of the original subdivision, the inductive procedure yields a subdivision in which all cones containing $X_\nu$ consist of this ray and a collection of rays added in the course of the subdivision. Note also that only vectors $v_i$ belonging to the valuation cone have been added in this procedure, and that in the construction we have labeled them by $G$-invariant divisors. Thus the set of rays labeled by colors remains unchanged.
\end{proof}

\begin{remark}
For toric varieties, or more generally horospherical varieties, the cone of $G$-invariant valuations coincides with the entire ambient space. For toroidal varieties, all cones of the fan are likewise contained in the valuation cone. Thus the condition concerning the intersection with the cone of $G$-invariant valuations is automatically satisfied, and one can simply pass to a sufficiently fine subdivision. For example, in the case of a simplicial fan, it suffices to apply barycentric subdivision twice in succession.
\end{remark}

Recall also the procedure of toroidalization.

\begin{definition}
A spherical variety is called \emph{toroidal} if no color contains a $G$-orbit in $X$.
\end{definition}

\begin{definition}
The \emph{toroidalization} of a spherical variety $X$ is the universal surjective $G$-equivariant morphism $f:\widetilde X\rightarrow X$, where $\widetilde X$ is a toroidal variety dominating $X$. That is, every dominant map from a toroidal variety to $X$ factors through $f$.
\end{definition}

The toroidalization of $X$ is determined by the colored fan whose colored cones have no colors and in which every cone $\mathcal C_{X,Y}$ of the fan of $X$ is replaced by its intersection $\mathcal C_{X,Y}\cap\mathcal V_X$ with the valuation cone $\mathcal V_X$. If a colored ray of $\mathcal C_{X,Y}$ was labeled by a color whose valuation image lies in $\mathcal V_X$, the ray is labeled by a $G$-invariant valuation instead.

\section{Description of the Movable Cone}

\begin{definition}
For $\nu\in\RGX{X}$, let $\mathcal P_{X,\widehat{\nu}}\subset\operatorname{Coker} A^*$ denote the monoid generated by the images under $\pi$ of the rays $X_\mu$, where $\mu\in\RGX{X}\setminus\{\nu\}$, and of all the rays $D\in\mathcal D$, with nonnegative integral coefficients.
\end{definition}

Thus, for $\nu\in\RGX{X}$,
\[
\mathcal P_{X,\widehat{\nu}}=\sum_{\mu\in\RGX{X}\setminus\{\nu\}}n_\mu X_\mu+\sum_{D\in\mathcal D}n_D D,\qquad n_\mu,n_D\in\mathbb Z_{\geq 0},
\]
where the sums are understood in the divisor class group.

\begin{theorem}\label{Mov}
Let $X$ be a spherical variety. Then the monoid $\MP_{\Mov(X)}$ of classes of effective Weil divisors whose complete linear systems have base loci of codimension at least $2$ is given by the formula
\[
\MP_{\Mov(X)}:=\MP_{\Eff}\cap\bigcap_{\nu\in\RGX{X}}\MP_{X,\widehat{\nu}},
\]
where $\nu\in\RGX{X}$ runs over all valuations corresponding to $G$-invariant prime divisors on $X$.
\end{theorem}

\begin{proof}
Let $E$ be a movable divisor as in the statement of the theorem, that is, an effective divisor whose base locus has codimension at least $2$. As before, after passing to a covering of $G$, we may assume that the sheaf $\mathcal O(E)$ is $G$-linearized. Then the base locus $\Bs(|E|)$ of its complete linear system is $G$-invariant. The condition that this locus has codimension at least $2$ is equivalent to the condition that $\Bs(|E|)$ contains none of the $G$-invariant divisors $X_\nu$, where $\nu\in\RGX{X}$.
Consequently, the restriction of the linear system to the dense $G$-orbit in the divisor $X_\nu$ is base-point-free, which is equivalent to $\pi(E)\in\mathcal P_{X,\widehat{\nu}}$. Indeed, the subspace of sections whose divisors contain $X_\nu$ is a proper $G$-invariant subspace. By complete reducibility, there exists a $B$-semi-invariant section in its $G$-invariant complement. The corresponding effective $B$-invariant divisor does not contain $X_\nu$. The converse is immediate.
Combining these conditions for all divisors $X_\nu$, where $\nu\in\RGX{X}$, and taking into account that $E$ is effective, we obtain the desired formula.
\end{proof}

Let $\mathcal P^{\mathbb Q}_{X,\widehat{\nu}}$ denote the rational cone generated by the monoid $\mathcal P_{X,\widehat{\nu}}$.

\begin{corollary}
Let $X$ be a spherical variety. Then the movable cone is given by the formula
\[
\Mov(X)=\Eff(X)\cap\bigcap_{\nu\in\RGX{X}}\mathcal P^{\mathbb Q}_{X,\widehat{\nu}},
\]
where the intersection is taken in $\Cl(X)_{\mathbb Q}$ and $\nu$ runs over all valuations corresponding to $G$-invariant prime divisors on $X$.
In particular, $\Mov(X)$ is a rational polyhedral cone.
\end{corollary}

\section{Decomposition of $\Mov$ into Chambers}

We now consider a decomposition of the movable cone into convex cones.
To this end, consider the set of quasiprojective spherical varieties $X_i$ that are isomorphic to $X$ outside subsets of codimension at least $2$.
It follows directly from the definition that the colored fan $\Sigma_i$ corresponding to $X_i$ has the same set of colored rays $\RG$, corresponding to $G$-invariant divisors, as the fan $\Sigma$ of $X$.
Note that, by Theorem~\ref{Mov}, we have $\Mov(X)=\Mov(X_i)$ for every $i$.

\begin{theorem}\label{dec}
Let $X$ be a projective spherical variety. The cone $\Mov(X)$ decomposes into cones $\Nef(\widetilde{X}_i)$ of nef divisors, where $\widetilde{X}_i$ are spherical varieties $G$-birational to $X$. This decomposition forms a fan, that is, any two cones $\Nef(\widetilde{X}_i)$ of the decomposition intersect along a common face, including on the boundary of $\Mov(X)$.
The cones of maximal dimension in the decomposition correspond to varieties $X_i$ that are isomorphic to $X$ outside subsets of codimension at least $2$.
\end{theorem}

We postpone the proof of the theorem in order to make some preliminary observations.
The following proposition, well known in the theory of toric varieties~\cite{Arzh}, gives a sufficient condition for the separatedness of a colored fan. For completeness, we include its proof.

\begin{proposition}\label{separation}
The convex cones $\mathcal C_{X,Y_1}$ and $\mathcal C_{X,Y_2}$ intersect along a common face, and the sets of indexed generators of the two cones lying in this face coincide, if and only if the relative interiors of the Gale dual cones $\mathcal P^{\mathbb Q}_{X,\widehat{Y}_1}$ and $\mathcal P^{\mathbb Q}_{X,\widehat{Y}_2}$ have nonempty intersection.
\end{proposition}

\begin{remark}
Note that the condition of the proposition is equivalent to the existence of an ample line bundle on the possibly nonseparated variety obtained by gluing the simple varieties $X_{Y_1}$ and $X_{Y_2}$. This yields a locally closed embedding of the variety into projective space and automatically implies its separatedness.
\end{remark}

\begin{proof}
Let $v_0 \in (\mathcal P^{\mathbb Q}_{X,\widehat{Y}_1} )^\circ \cap (\mathcal P_{X,\widehat{Y}_2}^{\mathbb Q} )^\circ$ be a vector.
Index the rays of the colored fan as follows: $I:=\RGX{X} \sqcup \MD$, $I_1:=\RGX{X,Y_1} \sqcup \mathcal D_{Y_1}$, and $I_2:=\RGX{X,Y_2} \sqcup \mathcal D_{Y_2}$. Denote the basis elements of the group $\mathbb Z\langle X_\nu, D\rangle$ corresponding to indices $\alpha \in I$ by $e_\alpha$.
On the one hand, we can choose a lift (not necessarily unique) of the class $v_0$ to the space $\mathbb Q\langle X_\nu, D\rangle$ of the form
$v=\sum_{\alpha \in I \setminus I_1} a_\alpha e^*_\alpha$, where $a_\alpha > 0$.

On the other hand, there exists $\chi \in \Xi(X)_{\mathbb Q}$ such that the element $v+A^*\chi$, representing the class $v_0$, can be written as
$v+A^*\chi=\sum_{\beta \in I \setminus I _2} b_\beta e^*_\beta$, where $b_\beta > 0$. Consider the vector
\[
A^*\chi=\sum_{\beta \in I \setminus I _2} b_\beta e^*_\beta
-\sum_{\alpha \in I \setminus I_1} a_\alpha e^*_\alpha,
\]
which can be used to separate the corresponding colored cones.
Let $\nu_0\in \mathcal C_{X,Y_2}$. Then $\nu_0 =Ae_0$, where $e_0=\sum_{\gamma \in I_2} c_\gamma e_\gamma$ and $c_\gamma\geq0$. We have
\[
\langle \nu_0, \chi \rangle = \langle A e_0, \chi \rangle = \langle e_0, A^*\chi \rangle.
\]
Moreover,
\[
\langle e_0, A^*\chi \rangle= \Big\langle \sum_{\gamma \in I_2} c_\gamma e_\gamma,
\sum_{\beta \in I \setminus I _2} b_\beta e^*_\beta  -\sum_{\alpha \in I \setminus I_1}
 a_\alpha e^*_\alpha \Big\rangle
\]
\[
= - \sum_{ \gamma \in I_2 \setminus I_1} c_\gamma a_\gamma \leq 0.
\]
Similarly, the linear form $\chi$ is nonnegative on $\mathcal C_{X,Y_1}$. The strict positivity of the coefficients implies that
\[
\mathcal C_{X,Y_1}\cap\chi^\perp
=\mathcal C_{X,Y_2}\cap\chi^\perp
=\langle Ae_\alpha:\alpha\in I_1\cap I_2\rangle_{\geqslant 0}.
\]
Thus the intersection of the cones is a common face, and the sets of indexed generators lying in this face coincide.

It is easy to see that the above arguments can be reversed. Indeed, let $\chi$ be a separating element, that is, one satisfying $\langle Ae_\alpha, \chi \rangle \leq 0$ for $\alpha\in I_1$ and $\langle Ae_\alpha, \chi \rangle \geq 0$ for $\alpha \in I_2$, with equality precisely for $\alpha \in I_1 \cap I_2$. By construction,
\[
A^*\chi=  \sum_{\alpha \in I_2\setminus  I_1} a_\alpha e^*_\alpha
 -\sum_{\beta \in I_1 \setminus I_2} b_\beta e^*_\beta + \sum_{\gamma \in I \setminus
  I_1 \cup I_2} c_\gamma e^*_\gamma,
\]
where $a_\alpha$ and $b_\beta$ are strictly positive (being the absolute values of the corresponding pairings $\langle \nu, \chi \rangle$), whereas the $c_\gamma$ may have either sign.
Consider the vector
\[
v=  \sum_{\alpha \in I_2 \setminus  I_1} a_\alpha e^*_\alpha + \sum_{\nu \in I \setminus
  I_1 \cup I_2} d_\nu e^*_\nu,
\]
whose class lies in $(\MP^{\mathbb Q}_{X,\widehat{Y}_1})^\circ$, where the $d_\nu$ are arbitrary strictly positive numbers satisfying $d_\nu\geq |c_\nu|+1$.
It is easy to see that the vector
\[
v-A^*\chi=  \sum_{\beta \in I_1 \setminus I_2} b_\beta e^*_\beta + \sum_{\nu \in I \setminus
  I_1 \cup I_2} (d_\nu-c_\nu) e^*_\nu
\]
represents the same class in $\Coker A^*$, and that this class belongs to $(\MP^{\mathbb Q}_{X,\widehat{Y}_2})^\circ$.
\end{proof}

\begin{definition}
A collection of colored cones that need not satisfy condition~\eqref{eq:separ} is called a \emph{colored quasi-fan}.
\end{definition}

The following proposition follows directly from Proposition~\ref{separation} and the separatedness condition for a colored fan.

\begin{proposition}
Consider a colored quasi-fan $\widetilde{\Sigma}$ for the same spherical homogeneous space as $X$, closed under taking colored faces. Suppose that the set of colored rays of $\widetilde{\Sigma}$ corresponding to $G$-invariant prime divisors coincides exactly with the set of rays generated by the vectors $\nu\in\RGX{X}$. The sets of colors of individual cones are subsets of the original set $\mathcal D$; the distribution of colors among the cones need not coincide with their distribution in $\Sigma$. Suppose that the cone $\mathcal P^{\mathbb Q}_{\widetilde{\Sigma}}$ of Theorem~\ref{Nef}, constructed from the quasi-fan $\widetilde{\Sigma}$, has nonempty interior in $\Cl(X)_{\mathbb Q}$ (is full-dimensional). Then $\widetilde{\Sigma}$ is a colored fan and corresponds to a spherical variety that is $G$-equivariantly isomorphic to $X$ on open subsets whose complements have codimension at least $2$. If, in addition, the support of $\widetilde{\Sigma}$ coincides with $\mathcal V_X$, then this variety is complete.
\end{proposition}

\begin{corollary}
The movable cone decomposes into chambers for which the equivalence class determines a model of a spherical variety isomorphic to $X$ outside a subset of codimension at least $2$.
\end{corollary}

\section{Big Divisors and Birational Models}

Recall the following definition.

\begin{definition}
A Weil divisor is called \emph{big} if the linear system of a sufficiently large multiple of it defines a birational map onto its image.
\end{definition}

The following statement is also well known.

\begin{proposition}
A Weil divisor on a normal projective variety is big if and only if a positive multiple of it is linearly equivalent to the sum of an ample Cartier divisor and an effective Weil divisor.
\end{proposition}

In particular, this criterion implies that every divisor in the interior $\Eff(X)^\circ$ of the rational cone of effective divisors is big. (Indeed, from any effective rational divisor in $\Eff(X)^\circ$ one can subtract an ample divisor $tA$ for a small positive $t\in\mathbb Q$ so that the resulting divisor $E$ remains in $\Eff(X)^\circ$. Hence $E+tA$ is big.) This also implies the convexity of the cone of big divisors.
If $D\in\Eff(X)^\circ$, then the image $Y$ of the map defined by the complete linear system of a sufficiently large and divisible multiple of $D$ has the same dimension as $X$ (and thus $\Lambda(X)=\Lambda(Y)$). Suppose now that $D\in\Mov(X)$.
Consider the subset $U\subset X$ consisting of the open $G$-orbit in $X$ and the open $G$-orbits in the $G$-invariant divisors on $X$. Then $H^0(U,\mathcal O_U(nD))=H^0(X,\mathcal O_X(nD))$, and $|nD|$ is base-point-free on $U$ for sufficiently large and divisible $n$.
Thus $Y$ can be recovered as $\overline{f(U)}$.
Some divisors on $U$ may be contracted by $f$; consequently, the rays of the fan of the variety $f(U)$ form a subset of $\Sigma(1)$.

\begin{proposition} (\cite{BrionGeom}) Let $X$ be a projective spherical $G$-variety. The class of a subvariety $Y\subset X$ in the Chow group can be represented by an effective $B$-invariant cycle. Since $X$ has finitely many $B$-orbits, the components of the support of any $B$-invariant cycle belong to a finite set of $B$-stable subvarieties.
\end{proposition}

\begin{lemma}\label{lem:finitegeneration}
Let $X$ be a normal projective spherical variety and let $D$ be a Weil divisor on $X$. Then the section ring
\[
R(X,D)=\bigoplus_{m\geqslant 0}H^0(X,\mathcal O_X(mD))
\]
is finitely generated.
\end{lemma}

\begin{lemma}\label{lem:noextraction}
Let $X$ be a normal projective variety and let $D$ be a big Weil divisor on $X$ whose section ring
\[
R(X,D)=\bigoplus_{m\geq 0}H^0(X,\mathcal O_X(mD))
\]
is finitely generated. Then the birational map $X\dashrightarrow X'=\operatorname{Proj}R(X,D)$ does not extract prime divisors, i.e., every prime divisor on $X'$ defines the same divisorial valuation of the function field as some prime divisor on $X$.
\end{lemma}
\begin{proof}
Let $s\in H^0(X,\mathcal O_X(mD))$, $m>0$, be a nonzero section, and let $D_s$ be its divisor of zeros. For the standard affine open subset $X'_s\subset X'$ we have
\[
\Bbbk[X'_s]=R(X,D)_{(s)}
=\bigcup_{k\geq0}\frac{H^0(X,\mathcal O_X(kmD))}{s^k}
=\Gamma(X\setminus\operatorname{Supp}D_s,\mathcal O_X).
\]
The last equality follows from the fact that every regular function on the complement of the support of the effective Weil divisor $D_s$ has poles of bounded order along its components. By normality of $X$ we obtain
\[
\Bbbk[X'_s]=\bigcap_{E\not\subset\operatorname{Supp}D_s}\mathcal O_{X,\eta_E},
\]
where $E$ ranges over the prime divisors on $X$ and $\eta_E$ denotes their generic points. This intersection of discrete valuation rings has finite character: for every nonzero rational function, only finitely many of the corresponding valuations take nonzero values. By a standard property of Krull domains, every valuation corresponding to a prime divisor on $X'_s$ belongs to this family of valuations. Since the open subsets $X'_s$ cover $X'$, the assertion follows.
\end{proof}

\begin{lemma}[{\cite{Brion}, \S~5.3, Propositions~1, 2}]\label{lem:momentfan}
Let $X$ be a projective spherical variety and let $D$ be an ample $B$-invariant Cartier divisor. Denote by $I$ the index set of the $B$-invariant prime divisors $E_i$ on $X$ and write $D=\sum_{i\in I}a_iE_i$. Set
\[
\Delta_D=\{\chi\in\Xi(X)_{\mathbb Q}\mid \langle Ae_i,\chi\rangle+a_i\geq0,\ i\in I\}.
\]
For a nonempty face $F$ of the polytope $\Delta_D$, set
\[
I_F=\{i\in I\mid \langle Ae_i,\chi\rangle+a_i=0\ \text{for all }\chi\in F\}.
\]
Then the colored fan of $X$ consists of the pairs
\[
(\mathcal C_F,\mathcal D_F),\qquad
\mathcal C_F=\langle Ae_i:i\in I_F\rangle_{\geqslant 0},\quad
\mathcal D_F=I_F\cap\mathcal D,
\]
for which $\mathcal C_F^\circ\cap\mathcal V_X\ne\varnothing$; here the indices of colors are identified with the colors themselves.
\end{lemma}
\begin{proof}
After choosing a rational $B$-semi-invariant section, the moment polytope differs from $\Delta_D$ only by a translation. By \cite{Brion} (\S~5.3, Propositions~1, 2), the colored cones are the inward normal cones of the corresponding faces, and the colors are determined by the equalities holding throughout the face. The normal cone is generated by the normals of all active inequalities.
\end{proof}

\begin{proposition}\label{neib}
Let $D\in\Mov(X)\cap\Eff(X)^\circ$ be a movable rational Weil divisor on a projective spherical variety $X$, and let $H_0$ be an ample Cartier divisor on $X$.
Then, for every sufficiently small rational $t>0$, the divisor $D_t=D+tH_0$ is movable, and the complete linear system of a sufficiently large and divisible multiple of it defines an immersion at the generic point of every prime divisor on $X$.
Moreover, there exists a projective model $X'$ that is $G$-equivariantly birational to $X$ and isomorphic to it in codimension one, such that the strict transform of $D_t$ on $X'$ is an ample rational Cartier divisor.
\end{proposition}
\begin{proof}
Let $U\subset X$ be the open subset defined as the complement of the $G$-orbits of codimension at least $2$.
Suppose that $D\in\Mov(X)\cap\Eff(X)^\circ$. Then $D$ is big, and the complete linear system of a sufficiently large multiple of it defines a birational map onto its image.
Replacing $D$, if necessary, by a positive integral multiple, we use the construction of Proposition~\ref{prop:movdef}, where $\dim Y=\dim X$, and assume that $D=\psi_*\widetilde f^*H$ and $\operatorname{codim}\Bs(|D|)\geq 2$. This replacement entails no loss of generality, since the assertion is invariant under positive scaling.
By assumption, $\Bs(|D|)$ does not meet $U$. Moreover, $U$ is smooth, since a normal variety is regular in codimension one and its singular locus is $G$-invariant. Hence $\psi$ is an isomorphism over $U$.

Consider the rational divisor $D_t=D+tH_0$ for a small rational $t>0$.
On $\widetilde X$ consider the divisor
\[
\widetilde D_t=\widetilde f^*H+t\psi^*H_0.
\]
Since both summands are semiample, the linear system of a sufficiently large and divisible multiple $m\widetilde D_t$ is base-point-free. Choose $m$ so that $mH$ and $mtH_0$ are very ample divisors.
Products of pullbacks of sections of these divisors define a linear subsystem of $|m\widetilde D_t|$ whose associated map is the composition of a map to a product of projective spaces and the Segre embedding. Since the second map to a projective space is the composition of $\psi$ with the embedding defined by $|mtH_0|$, this subsystem defines a locally closed immersion on $\psi^{-1}(U)$.
Since $\psi_*\widetilde D_t=D_t$ and $X$ is normal, global sections of the multiple $m\widetilde D_t$ define on $U$ a subsystem of the complete linear system $|mD_t|$. Consequently, the complete linear system $|mD_t|$ also defines a locally closed immersion on $U$, in particular at the generic point of every prime divisor on $X$.
Furthermore, $D_t$ is movable, since the sum of a movable divisor and an ample divisor is movable.

By Lemma~\ref{lem:finitegeneration}, the section ring $R(X,mD_t)$ is finitely generated. Set $X'=\operatorname{Proj}R(X,mD_t)$, increasing $m$ if necessary so that the section ring of the corresponding multiple is generated in degree one. The variety $X'$ is normal, and the map $X\dashrightarrow X'$ is birational and $G$-equivariant. By Lemma~\ref{lem:noextraction}, it does not extract divisors, and by the preceding argument it does not contract divisors on $X$. Consequently, it is an isomorphism in codimension one. The strict transform of $mD_t$ on $X'$ is linearly equivalent to the ample divisor corresponding to $\mathcal O_{X'}(1)$: this follows from the coincidence of the spaces of sections and the absence of fixed divisorial components on $X$.

In particular, choosing a $B$-invariant representative $D_t=\sum_{i\in I}a_i(t)E_i$, we can describe the colored fan of $X'$ directly in terms of the original divisorial data of $X$. To this end, we construct the polytope $\Delta_{D_t}$ of Lemma~\ref{lem:momentfan} with coefficients $a_i(t)$ and determine the set of active indices $I_F$ for each of its nonempty faces $F$. The colored fan of $X'$ consists of the pairs $(\langle Ae_i:i\in I_F\rangle_{\geqslant 0},I_F\cap\mathcal D)$ for which the relative interiors of the cones meet $\mathcal V_X$. Here we have used the fact that all the original divisorial valuations and colors are preserved under an isomorphism in codimension one; passing to a multiple merely rescales the polytope.

Letting $t$ tend to zero, we conclude that every neighborhood of the class of $D$ contains the class of a movable divisor such that the linear system of a sufficiently large and divisible multiple of it defines an immersion at the generic point of every prime divisor on $X$.
\end{proof}

\begin{remark}
We relate the description of the colored fan of the model $X'$ in Proposition~\ref{neib} to Gale duality, using the notation of Proposition~\ref{separation}. Set $I:=\RGX{X}\sqcup\MD$ and write $D_t=\sum_{\alpha\in I}a_\alpha e^*_\alpha$, where $e^*_\alpha$ are the basis elements of the group of $B$-invariant divisors dual to $e_\alpha$.
For a subset $J\subseteq I$, the existence of a nonempty face $F$ of the polytope $\Delta_{D_t}$ with active index set $I_F=J$ is equivalent to the existence of $\chi\in\Xi(X)_{\mathbb Q}$ such that
\[
D_t+A^*\chi=\sum_{\alpha\in I\setminus J}c_\alpha e^*_\alpha,\qquad c_\alpha>0.
\]
Indeed, such a point $\chi$ can be chosen in the relative interior of $F$, and conversely, such a point determines a face with the specified active index set. By the definition of the map $\pi$ (after tensoring with $\mathbb Q$), this is equivalent to the condition
\[
\pi(D_t)\in\bigl\langle \pi(e^*_\alpha):\alpha\in I\setminus J\bigr\rangle_{>0}.
\]
If, in addition, $\langle Ae_\alpha:\alpha\in J\rangle_{>0}\cap\mathcal V_X\ne\varnothing$, then by Lemma~\ref{lem:momentfan} this face corresponds to a $G$-orbit $Y\subset X'$, and $J=\RGX{X',Y}\sqcup\mathcal D_Y$. Identifying the divisorial data of $X$ and $X'$, the last condition on the class can be written in the notation of Proposition~\ref{separation} as
\[
\pi(D_t)\in\bigl(\MP^{\mathbb Q}_{X',\widehat{Y}}\bigr)^\circ.
\]
Thus the colored cones of the model and their colors can be recovered from membership of the divisor class in the relative interiors of the Gale dual cones, without constructing the polytope itself.
\end{remark}

\section{Comparison of Birational Models and Their Nef Cones}

\begin{lemma}\label{lem:iso}
Let $X_1$ and $X_2$ be projective varieties that are $G$-equivariantly isomorphic on open subsets $U_1$ and $U_2$ whose complements have codimension at least two. Let $\mathcal L_1,\mathcal L_2$ be ample line bundles on $X_1$ and $X_2$, respectively, representing the same element under the identification of the class groups induced by the inclusions $U_i\subset X_i$ and the isomorphism $f:U_1\rightarrow U_2$. Then $X_1$ and $X_2$ are isomorphic.
\end{lemma}
\begin{proof}
Replacing $\mathcal L_i$ by their powers $\mathcal L_i^l$, we may assume that both line bundles are very ample. Let $U_i\subset X_i$ be the subsets on which the isomorphism $f:U_1\rightarrow U_2$ is defined. We then have identifications
\[
H^0(X_1,\mathcal L_1)=H^0(U_1,\mathcal L_1)=H^0(U_2,\mathcal L_2)=H^0(X_2,\mathcal L_2),
\]
where the first and last maps are isomorphisms because the complements of $U_i$ in $X_i$ have codimension at least two, and the middle map is an isomorphism because the elements of the class groups coincide. Since $\mathcal L_i$ are very ample, the maps $g_i:X_i\rightarrow\mathbb P(H^0(X_i,\mathcal L_i)^*)$ are embeddings. The coincidence of the linear systems implies that $g_i$ agree on the subsets $U_i$. Thus $g_1(X_1)=\overline{g_1(U_1)}=\overline{g_2(U_2)}=g_2(X_2)$, yielding an isomorphism between the varieties $X_i$.
\end{proof}

We shall need the following lemma in our study of linear systems.
\begin{lemma}\label{fiber} (\cite{Debarre})
Let $X$ be a complete variety and let $f_i:X\rightarrow Y_i$, $i=1,2$, be surjective morphisms such that every fiber of $f_1$ is contained in a fiber of $f_2$; that is, for every $y_1\in Y_1$ there exists $y_2\in Y_2$ such that $f_1^{-1}(y_1)\subset f_2^{-1}(y_2)$. Suppose that $Y_1$ is normal. Then there exists a morphism $g:Y_1\rightarrow Y_2$ such that $f_2=g\circ f_1$. If the fibers of $f_1$ and $f_2$ coincide and $Y_2$ is normal, then $g$ is an isomorphism.
\end{lemma}
\begin{proof}
Consider the map $\Phi=f_1\times f_2:X\rightarrow Y_1\times Y_2$ and denote its image by $\Gamma$ (it is closed, being the image of a proper variety, and projects surjectively onto each $Y_i$). The containment of the fibers of $f_1$ in those of $f_2$ implies that the projection from $\Gamma$ onto $Y_1$ is bijective. Since surjectivity has already been established, it suffices to prove injectivity. Indeed, suppose two points of $\Gamma$ map to $y_1$. They can be written as $\Phi(x)$ and $\Phi(\tilde x)$, with $f_1(x)=f_1(\tilde x)$. By the hypothesis on the fibers, $f_2(x)=f_2(\tilde x)$, proving that $\Phi(x)=\Phi(\tilde x)$. Since $Y_1$ is normal and the field $k$ has characteristic zero, the restriction of the projection $p_1$ to $\Gamma$ is an isomorphism onto $Y_1$. We may take $g=p_2\circ(p_1|_\Gamma)^{-1}:Y_1\rightarrow Y_2$.

If the fibers of $f_1$ and $f_2$ coincide and $Y_2$ is normal, the inverse of $g$ is $p_1\circ(p_2|_\Gamma)^{-1}$, since $p_2|_\Gamma$ is also an isomorphism.
\end{proof}

To compare nef cones, we use the following lemma.
\begin{lemma}\label{comm}
Let $p:\widetilde X\rightarrow X$ be a surjective morphism of normal projective varieties with connected fibers. Suppose that $\Nef(\widetilde X)$ is polyhedral and sufficiently divisible positive multiples of its elements are base-point-free. Let $\widetilde D$ belong to the smallest face of $\Nef(\widetilde X)$ containing $p^*\Nef(X)$. Let $Y$ be the image of $\widetilde X$ under the map $\varphi_{\widetilde D}$ defined by the complete linear system $|\widetilde D|$. Then there exists a map $g:X\rightarrow Y$ such that $\varphi_{\widetilde D}=g\circ p$. Moreover, $\widetilde D=p^*g^*\mathcal O_Y(1)$ is the pullback of the very ample divisor $\mathcal O_Y(1)$ on $Y$ (and hence the pullback of the nef divisor $g^*\mathcal O_Y(1)$ on $X$).
\end{lemma}
\begin{proof}
Let $D$ be a very ample divisor on $X$. Then $p^*D$ belongs to the relative interior of $p^*\Nef(X)$. Let $F$ be the smallest face of $\Nef(\widetilde X)$ containing $p^*\Nef(X)$.

Let $C$ be a curve whose intersection number with $p^*D$ is zero. Since $F$ is a face of the nef cone, which is dual to the cone of effective curves, all elements of $F$ have zero intersection with $C$. As $\widetilde D\in F$, the intersection number of $C$ with $\widetilde D$ is also zero. The curves having zero intersection with $p^*D$ are precisely those contained in fibers of $p$. The fiber $p^{-1}(p(x))$, being connected, can be recovered as the set of points of $\widetilde X$ that can be joined to $x$ by a chain of curves having zero intersection with $p^*D$.

Replace $\widetilde D$ by a sufficiently divisible multiple so that the complete linear system defining $\varphi_{\widetilde D}$ is base-point-free. Thus $\varphi_{\widetilde D}$ is regular on all of $\widetilde X$, and completeness implies that its image $Y$ is closed in projective space.

In turn, the curves having zero intersection with $\widetilde D$ are precisely those contained in fibers of $\varphi_{\widetilde D}$. Two points of $\widetilde X$ lie in the same connected component of a fiber of $p$ (respectively, $\varphi_{\widetilde D}$) if and only if they can be joined by a chain of irreducible curves, each having zero intersection with $p^*D$ (respectively, $\widetilde D$). Since the fibers of $p$ are connected, the fibers of $\varphi_{\widetilde D}$ contain those of $p$. Thus the hypotheses of Lemma~\ref{fiber} hold for $f_1=p$ and $f_2=\varphi_{\widetilde D}$, and the desired map $g$ exists.

\[
\xymatrix{
 & & \widetilde X \ar[dl]_(0.5){\varphi_D} \ar[d]_p \ar[dr]^(0.5){\varphi_{\widetilde D}} & \\
 & \mathbb P(p^*H^0(X,\mathcal O(D))^*) & \ar@{_{(}->}[l]\ar[r]^(0.25){g} X & Y\rightarrow\mathbb P(H^0(X,\mathcal O(\widetilde D))^*) & ,}
\]
By construction of $\varphi_{\widetilde D}$, we have $\widetilde D=\varphi_{\widetilde D}^*\mathcal O_{\mathbb P^l}(1)$, where $\mathbb P^l=\mathbb P(H^0(X,\mathcal O(\widetilde D))^*)$. Thus $\widetilde D$ is the pullback of a very ample divisor on $Y$, and also the pullback under $p$ of the nef divisor $g^*\mathcal O_Y(1)$ on $X$.
\end{proof}

\begin{proposition}\label{face}
Let $p:\widetilde X\rightarrow X$ be a surjective morphism of projective varieties with connected fibers. Suppose that $\Nef(\widetilde X)$ is polyhedral and sufficiently divisible positive multiples of its elements are base-point-free. Then $p^*\Nef(X)$ is a face of $\Nef(\widetilde X)$.
\end{proposition}
\begin{proof}
Choose a class $D$ in the relative interior of $\Nef(X)$. The map $p^*$ embeds $\Nef(X)$ into $\Nef(\widetilde X)$. Let $F$ be the smallest face of $\Nef(\widetilde X)$ containing $p^*D$. We shall show that $F=p^*\Nef(X)$.

For any $\widetilde D=p^*D\in(p^*\Nef(X))^\circ$, the intersection number with every curve in a fiber of $p$ is zero, so $\widetilde D$ lies on the boundary of $\Nef(\widetilde X)$; since $p^*\Nef(X)$ is a convex cone, it is contained in the face $F$. On the other hand, for every curve $C\subset\widetilde X$ whose image $p(C)$ is not a point, the intersection number $(C,p^*D)=(p(C),D)$ is strictly positive, since $D$ is ample. Consider the linear subsystem $p^*|D|=p^*H^0(X,\mathcal O(D))$ of $H^0(X,\mathcal O(\widetilde D))$, and let $\varphi_{p^*|D|}:\widetilde X\rightarrow\mathbb P(p^*H^0(X,\mathcal O(D))^*)$ be the map defined by this subsystem. Denote by $\widehat\varphi_{p^*|D|}$ the map with connected fibers obtained by Stein factorization. Since the elements of this system have zero intersection with curves in the fibers of $p$, and the fibers themselves are connected, the system does not separate points in the connected components of the fibers. By Lemma~\ref{fiber}, $\widehat\varphi_{p^*|D|}$ factors through $p$. Since $|D|$ defines an embedding of $X$ into projective space, the image of $\widehat\varphi_{p^*|D|}$ coincides with $X$.

We now show that every divisor $\widetilde D\in F^\circ$ belongs to $(p^*\Nef(X))^\circ$. Recall that, by Kleiman's criterion, the polyhedral cone $\Nef(\widetilde X)$ is dual to the cone of effective curves with respect to the intersection pairing. The face $F$ is dual to $F^\bot$, the cone of effective curves having zero intersection with every element of $F$ (equivalently, with an element of $F^\circ$).

Apply Lemma~\ref{comm} to $\widetilde D$ in the smallest face $F$ containing $p^*\Nef(X)$. We obtain $\widetilde D=p^*g^*\mathcal O_Y(1)$, so $\widetilde D\in p^*\Nef(X)$, as required.
\end{proof}

\begin{proposition}\label{lem:nef}
Let $X_1$ and $X_2$ be projective spherical varieties that are $G$-equivariantly isomorphic on open subsets $U_1$ and $U_2$ whose complements have codimension at least two. Identify the class groups of $X_i$ via the inclusions $j_i:U_i\rightarrow X_i$. Then the intersection of $\Nef(X_1)$ and $\Nef(X_2)$ in $\Cl(U_1)_{\mathbb Q}\cong\Cl(U_2)_{\mathbb Q}$ is a common face of these cones.
\end{proposition}
\begin{proof}
As before, we have an identification $f:U_1\rightarrow U_2$ of the open subsets. Pass to a common resolution $\widetilde X$ of the indeterminacy of $f$, viewed as a birational map between the $X_i$, and denote the corresponding morphisms by $p_i:\widetilde X\rightarrow X_i$.

Let $D_i\in\Nef(X_i)$ be rational classes of nef divisors whose restrictions to $U_1\cong U_2$ coincide. Passing to a common sufficiently divisible positive multiple, we may assume that $D_i$ are Cartier divisors whose complete linear systems are base-point-free. By normality of $X_i$ and the codimension of the complements of $U_i$, the spaces of global sections of these divisors are identified by restriction to $U_i$. Hence the corresponding morphisms have the same image $Y$ and agree on the common open subset $U$. Bigness of the divisors is not required here.

Set $F_i=p_i^*\Nef(X_i)$. As shown in Proposition~\ref{face}, this is a face of $\Nef(\widetilde X)$. Pass to a common sufficiently large multiple so that the complete linear systems of $p_1^*D_1$ and $p_2^*D_2$ are base-point-free on $\widetilde X$. The maps $\varphi_{p_i^*|D_i|}$ have the same image $Y=\overline{\varphi_{|D_1|}(U)}=\overline{\varphi_{|D_2|}(U)}$. Since the maps $\varphi_{|D_1|}p_1$ and $\varphi_{|D_2|}p_2$ agree on $U$, they agree on all of $\widetilde X$. On $\widetilde X$ they are identified with $\varphi_{p_1^*|D_1|}$ and $\varphi_{p_2^*|D_2|}$, respectively.

By construction, $D_i=\varphi^*_{|D_i|}\mathcal O_Y(1)$, where $\mathcal O_Y(1)$ is a very ample line bundle on $Y$. Since $\varphi_{|D_1|}\circ p_1=\varphi_{|D_2|}\circ p_2$, we have $p_1^*D_1=p_2^*D_2$.

We have the following commutative diagram.
\[
\xymatrix@R=3.0pc@C=1.5pc{
 & \widetilde{X} \ar[dl]_{p_1} \ar[dr]^{p_2}
   \ar@/^1.8pc/[dd]^{\raisebox{-17pt}{$\scriptstyle\widehat{\varphi}_{|\widetilde{D}|}$}}
   \ar[ddl]^(0.65){\raisebox{-18pt}{$\scriptstyle\varphi_{p_1^*|D_1|}\mkern-12mu$}}
   \ar[ddr]_(0.65){\raisebox{-18pt}{$\scriptstyle\mkern-12mu\varphi_{p_2^*|D_2|}$}} & \\
 X_1 \ar[d]_{\varphi_{D_1}\mkern-4mu} &
 U \ar@{_{(}->}[l]_{j_1} \ar@{^{(}->}[r]^{j_2}
   \ar@{_{(}->}[u]_{\widetilde{j}} &
 X_2 \ar[d]^{\mkern-4mu\varphi_{D_2}} \\
 \mathbb P(H^0(\widetilde{X},\mathcal O(p_1^*D_1))^*) &
 Y \ar@{_{(}->}[l] \ar@{^{(}->}[r] &
 \mathbb P(H^0(\widetilde{X},\mathcal O(p_2^*D_2))^*)
}
\]

Thus we have shown that if two divisor classes $D_1$ and $D_2$ coincide on $U$, then the classes $p_1^*D_1$ and $p_2^*D_2$ coincide on $\widetilde X$ and belong to the intersection of the faces $F_1:=p_1^*\Nef(X_1)$ and $F_2:=p_2^\ast\Nef(X_2)$ of $\Nef(\widetilde X)$. Hence their intersection $F=F_1\cap F_2$ is also a face of $\Nef(\widetilde X)$ and is a face of each cone $F_i$. Note also that if $\widetilde D\in F$, there exist divisors $D_i$ on $X_i$ such that $p_1^*D_1=p_2^*D_2$, and on the common open subset $U$ of $\widetilde X$, $X_1$, and $X_2$, the classes of $D_i$ coincide with the restriction of the class of the effective divisor $\widetilde D$. Recall that $(p_i)_*$ define embeddings of $p_i^*\Nef(X_i)$ into $\Cl(X_i)\cong\Cl(U)$. This means that the class $\widetilde D\in F$ can be represented compatibly under the embeddings
\[
F\subset p_i^*\Nef(X_i)\xrightarrow{\quad(j_i)_*(p_i)_*\quad}\Cl(U).
\]
This shows that the images of the cones $p_i^*\Nef(X_i)$ in $\Cl(U)$ (which coincide with the images of $\Nef(X_i)$, respectively) intersect precisely along the image of their common face $F$, as required.
\end{proof}

\begin{proof}[Proof of Theorem~\ref{dec}]
Consider the cones $\Nef(X_i)$ in $\Mov(X)$, where $X_i$ range over all models with the same set of colored rays $\RGX{X}$ marked by $G$-divisors. The models $X_i$ and $X$ contain a common $G$-invariant open subset consisting of the open $G$-orbit in $X$ and the dense orbits in the $G$-invariant divisors. There are clearly only finitely many such models, since they correspond to grouping the finitely many rays $\RGX{X}$ and colors into subsets defining the cones of the fan. By Lemma~\ref{lem:iso} and Proposition~\ref{lem:nef}, the interiors of distinct pointed full-dimensional cones $\Nef(X_i)$ are disjoint, and the cones themselves intersect along common faces.

By Proposition~\ref{neib}, the set of divisors in $\Mov(X)\cap\Eff(X)^\circ$ lying in one of the chambers $\Nef(X_i)$ is dense in $\Mov(X)$: every movable class can be approximated by big movable classes by adding small positive multiples of an ample divisor, and Proposition~\ref{neib} applies to these classes. Since the finite union of the cones $\Nef(X_i)$ is closed, it follows that
\[
\Mov(X)=\bigcup_{U\subset X_i}\Nef(X_i).
\]
\end{proof}

Thus Gale duality makes it possible to describe the movable cone of a spherical variety and its decomposition into chambers directly in terms of colored fans. The construction of equivariant birational models and the comparison of their nef cones thereby reduce to problems in convex geometry in the space of divisor classes.

\end{document}